\documentstyle[11pt]{article}
\newcommand{\Ibb}[1]{ {\rm I\ifmmode\mkern
            -3.6mu\else\kern -.2em\fi#1}}
\newcommand{\ibb}[1]{\leavevmode\hbox{\kern.3em\vrule
     height 1.2ex depth -.3ex width .2pt\kern-.3em\rm#1}}

\newcommand{\Rl}{{\Ibb R}}
\newcommand{\Nl}{{\Ibb N}}
\newcommand{\be}{\begin{eqnarray}}
\newcommand{\ee}{\end{eqnarray}}
\newcommand{\bez}{\begin{eqnarray*}}
\newcommand{\eez}{\end{eqnarray*}}
\renewcommand{\O}{\Omega}
\newcommand{\A}{{\cal A}}
\newcommand{\X}{{\cal X}}
\newcommand{\we}{\wedge}
\newcommand{\oa}{\otimes}
\newcommand{\bu}{\bullet}
\newcommand{\lb}{[\![}
\newcommand{\rb}{]\!]}
\newcommand{\pa}{\partial}
\newcommand{\na}{\nabla}
\newcommand{\To}[1]{\stackrel{#1}{\longrightarrow}}
\newcommand{\ins}{\leavevmode\hbox{\kern.pt\vrule height .2pt depth .pt 
    width 1.2ex \vrule height 1.2ex depth .pt width .2pt\kern2.2pt}}

\begin{document}
\begin{center}

{\Large\bf  NONCOMMUTATIVE GEOMETRY AND ITS RELATION TO\\
\vskip.1cm

STOCHASTIC CALCULUS AND SYMPLECTIC MECHANICS}
\vskip1cm

A. DIMAKIS$^\dagger$ and C.TZANAKIS$^\ddagger$\\
\vskip.2cm

$\dagger$ Institut F\"ur Theoretische Physik, Bunsenstr. 9, 37073 G\"ottingen, 
Germany\\

$\ddagger$ University of Crete, 74100 Rethymnon, Crete, Greece
\end{center}

\vskip1.cm

\centerline{\bf 1. INTRODUCTION}
\renewcommand{\theequation} {1.\arabic{equation}}
\setcounter{equation}{0}
\vspace*{.5cm}

\noindent
This paper rests on three, a priori quite distinct domains, namely, 
noncommutative geometry, stochastic calculus (StC) and symplectic 
mechanics (SyM). Hence we give below a very brief outline of some basic facts
concerning StC and SyM. 

StC is motivated by the desire to put in a firm mathematical basis, 
physically relevant but ill-defined differential equations (e.g.\ Langevin's
equation in Brownian motion), that finally have been interpreted as
stochastic differential equations. Accordingly, 
the appropriate tool for their study has 
been  developped, namely, stochastic integration of It\^o,
Stratonovich and others (e.g.\ \cite{Ar}, \cite{TG}, \cite{CG}).
The essential point to be emphasized in this context, is that stochastic
differentiation involves 2nd order derivatives, e.g.\ in It\^o's
formula (in its simplest form) for the differential of a function $f$
of a Wiener stochastic process $X_t$ (e.g.\ \cite{Ar} \S5.3, \cite{CG}
\S4.3)
\be 
df(X_t,t) = dt(\pa_t f+{\gamma\over2}\pa_x^2 f)+dX_t\,\pa_x f\;.\label{ito}
\ee
It is a general feature of StC that stochastic differentiation obeys a
generalized Leibniz rule
\be d(X_tY_t) = dX_t\,Y_t+dY_t\,X_t + d\lb X_t,Y_t \rb \label{sm} \ee
(e.g.\ \cite{Em} p.134), where the 3rd term is the so-called bracket of the 
processes (semimartingales), related to the quadratic variation of a process 
$X_t$
\be 
\lb X_t,X_t \rb =\lim_{k\to\infty}\sum_k(X_{t_k}-X_{t_{k-1}})^2 \label{1.3}
\ee
the limit being taken with respect to the width of the partition of the (time)
interval, where $t$ varies.

There have been efforts to extend these considerations to a semimartingale 
theory on differentiable manifolds, by appropriately introducing basic 
concepts like vector field, connection etc.\ (\cite{Em}, \cite{Me}).
We will come back to this later on.

In SyM we consider an even-dimensional manifold $M$ extended by $\Rl$, 
$(M\times\Rl,\omega)$, where $\omega$ is a closed 2-form of maximal
rank
\be  \omega={1\over2}dx_i\we dx_j\,\omega_{ij} + dt\we dH\;,\label{1.4} \ee
that is $\omega_{ij}$ is nondegenerate and $H$ is a 
function\footnote{In physics $M\times\Rl$ is the extended phase-space of
a physical system and $H$ its Hamiltonian.}.
Then $\omega$ has a 1-dimensional kernel, which by definition gives the 
Hamiltonian vector fields $X_H$ using the insert operator $\ins$
\be   X_H:\qquad X_H\ins\omega=0\;. \label{1.5}  \ee
Then Hamilton's equations are given by the action of $X_H$ on a smooth
function $A$ defined on $M\times\Rl$
\be X_H(A)=0\qquad \iff \qquad 	\pa_t A=-\{H,A\}\;,\label{hd} \ee
$\{\,,\}$ being the Poisson bracket of two functions. This is the starting 
point of both Hamiltonian dynamics (e.g.\ \cite{VA} ch.9) and statistical
mechanics. In the latter, making approximations based on iteration schemes
applied to (\ref{hd}), kinetic equations are derived, giving the
(probabilistic) time evolution of physical quantities in a system with many 
degrees of freedom (e.g.\ \cite{Li} ch.II). These are often generalized 
diffusion (or Fokker-Planck) equations (e.g.\ \cite{GrTz}, \cite{DiTz}
\S2)
\be 
\pa_t A = -\{H_0,A\}+\pa_i(a^{ij}\pa_j A)-\beta a^{ij}(\pa_j H_0)
\pa_i A\;, \label{fp}
\ee
summation convention always assumed\footnote{The terminology stems from
the fact that in the case of Brownian motion, time evolution is given by the
superposition of a diffusion process and a systematic drift (friction) term
--- 2nd and 3rd term in (\ref{fp}) respectively (cf.\ \cite{CG} ch.I).}. 
Let us notice here that it is possible to {\em derive kinetic equations via
StC, however at the expence of modifying Hamiltonian dynamics} (\cite{DiTz}
Section 3, \cite{CG} \S4.3). In fact {\em a common characteristic} of both
statistical mechanics and StC {\em is the appearence of  2nd order 
differential operators that cannot be given easily a geometrical meaning} in 
a natural way. Therefore, {\em our objective is to formulate an appropriate
geometric framework for this to be possible and develop SyM in it by retaining 
the Hamiltonian scheme and making only a minimal modification of the 
differential structure of the manifold}. Specifically we shall retain the 
manifold structure as this is encoded in the ordinary algebra of smooth 
functions on it and modify the differential calculus (DC). This leads us to 
noncommutative geometry and we shall see that tensor calculus and SyM can be 
developped in this framework as a direct extension of what is known in 
ordinary differential geometry.
\vspace*{1cm}

\centerline{\bf 2. NONCOMMUTATIVE DIFFERENTIAL CALCULUS ON 
MANIFOLDS\footnote{Lack of space does not allow for detailed proofs, which
will be presented elsewhere.}}
\renewcommand{\theequation} {2.\arabic{equation}}
\setcounter{equation}{0}
\vspace*{.5cm}

\noindent
We consider an $(N+1)$-dimensional manifold $M\times\Rl$ and the ordinary 
(commutative) algebra of (smooth) functions on it, $\A$. The universal
differential enveloppe of $\A$ is an $\Nl_0$-graded algebra $\O_u=
\bigoplus_{k\in\Nl_0}\O^k_u$ with $\O^0_u=\A$ equipped with an exterior
derivative operator $d_u$, satisfying well-known axioms (e.g.\ \cite{BDMH})
\begin{description}
\item[-] $d_u1=0$, $d_u^2=0$.
\item[-] $d_u(\psi\psi')=(d_u\psi)\psi'+(-1)^k\psi(d_u\psi')$, $\psi\in\O^k_u$.
\item[-] $\A$ and $d_u\A$ generate $\O_u$.
\end{description}
A simple representation of forms $\phi\in\O^r_u$, as functions on $(M\times
\Rl)^{r+1}$ is given by (\cite{DiTz})
\bez
(\phi\psi)(x_0,\ldots,x_{r+s})=\phi(x_0,\ldots,x_r)\psi(x_{r+1},\ldots,
x_{r+s})\;, 
\eez
\bez
(d_u\phi)(x_0,\ldots,x_{r+1}) = \sum_{k=0}^{r+1}(-1)^k\phi(x_0,\ldots,
x_{k-1},x_{k+1},\ldots,x_{r+1})\;,
\eez
$\psi\in\O^s_u$. Then we readily get $(d_uf)g\neq g\,d_uf$ in general, hence
the name noncommutative DC. Universality of $(\O_u,d_u)$
means that for any other DC $(\O,d)$ over $\A$, it can be proved that 
there exists
a unique graded algebra homomorphism
\be
\pi:\O_u\to\O\;,\qquad \pi|_\A=\mbox{id}_\A\;,\qquad d\circ\pi = 
\pi\circ d_u\label{2.1}
\ee
conserving the grade of forms. It is possible to make 1-forms $\O^1_u$ 
a commutative, associative algebra, by defining 
\be (\alpha\bu_u\beta)(x,y):= - \alpha(x,y)\,\beta(x,y)\;.\label{2.2} \ee
This implies 
\be
(d_uf)\bu_u \alpha =  [ f,\alpha ]:= f\alpha -\alpha f\;,\quad
f\in\A\;.
\ee
Thus {\em $\bu_u$ measures deviations from ordinary DC
due to noncommutativity}. In what follows we consider a {\em minimal}
deformation of it, namely a DC $(\O,d)$ on $M\times\Rl$
in which
\be
df\bu dg\bu dh=0\;,\qquad dt\bu df =0\;, \qquad f,g,h\in\A\;.\label{bu}
\ee
Such a DC can be obtained by taking the quotient of $\O_u$
by the differential ideal $J$, generated by
\bez 
d_uf\bu_u d_ug\bu_u d_uh\;,\qquad d_ut\bu_u d_uf \;, \qquad f,g,h\in\A\;,
\eez
that is 
\be
\pi:\O_u\to\O = \O_u/J\;,\quad\qquad \pi(\alpha)\bu\pi(\beta) = 
\pi(\alpha\bu_u \beta)\;, \label{pi}
\ee
$\pi$ being the canonical projection of $\O_u$ on $\O$. By (\ref{bu})
\be d(fg) = (df)g + (dg) f + df\bu dg \label{2.6}\ee
formally identical to (\ref{sm}) with $df\bu dg$ corresponding to
$d\lb f, g \rb$ (see also section 5 below).

$\O^1$ is an $\A$-bimodule. Then it is natural to define the space of 
vector fields $\X$ as its dual $\A$-bimodule\footnote{We define $\X$
as the left $\A$-module, dual to $\O^1$ taken as a right $\A$-module
only. The full $\A$-bimodule structure of $\X$ is then defined by
(\ref{bi}).}
\be 
\langle fXh,\alpha \rangle = f\langle X,h\alpha\rangle\;,\qquad
f,h\in\A\;,\;\;X\in\X\;,\;\;\alpha\in\O^1\;. \label{bi}
\ee
$\langle\,,\rangle$ denoting duality  and the action of $X$ on $\A$ is 
naturally given by
\be    X(f) :=  \langle X,df\rangle   \ee
implying
\be 
(fX)(g) = f\,X(g)\;,\qquad (Xf)(g)= X(fg)-X(f)\,g\;.  
\ee
It can be shown that $\X$ is characterized by $X\in\X\iff$
\be 
X(fgh)=fX(gh)+gX(fh)+hX(fg)-fgX(h)-fhX(g)-ghX(f)\label{xfgh}
\ee
a condition which is equivalent to 
\be  X(f^3)= 3 fX(f^2)-3f^2X(f)\;. \label{xt} \ee
Moreover (\ref{bu}) gives 
\be   X(tf)=tX(f)+fX(t)\;. \label{stx}\ee
Duality and (\ref{xfgh}), (\ref{xt}) give in local coordinates $(x^\mu,t)$
\be
X=X^t\pa_t + X^\mu \pa_\mu + {1\over2}X^{\mu\nu} \pa_\mu\pa_\nu\;,\label{xc}
\ee
\be 
df = dt\,\pa_tf + dx^\mu\pa_\mu f +{1\over2}\xi^{\mu\nu}\pa_\mu\pa_\nu f\;,
\qquad \xi^{\mu\nu}:=dx^\mu\bu dx^\nu\label{df}
\ee
in an obvious notation. Because of (\ref{xc}), (\ref{df}) we may call 
$(\O,d)$ a 2nd order DC\footnote{In StC on manifolds,
(\ref{xt}) has been used as a characterization of vector fields (\cite{Em}
Lemma 6.1).}. Moreover (\ref{df}) reduces formally to (\ref{ito}) in 
1-dimension, if $\xi^{\mu\nu}$ is proportional to $dt$ (see section 4 below).
Thus on $\X,\,\O^1$ we have respectively {\em noncovariant} bases
$\{\pa_t,\,\pa_\mu,\,\pa_\mu\pa_\nu\}$, $\{dt,\,dx^\mu,\,\xi^{\mu\nu}\}$.
To get {\em covariant} expressions, the concept of a connection has to be 
introduced.
\vspace*{1cm}

\centerline{\bf 3. CONNECTIONS}
\renewcommand{\theequation} {3.\arabic{equation}}
\setcounter{equation}{0}
\vspace*{.5cm}

\noindent
By considering the middle $\A$-linear tensor product $\oa$ of $\O^1$ with any
$\A$-bimodule $M$ (see e.g.\ \cite{Ja} \S3.7) it is possible to define a
(right) connection and on its dual $M^\ast$, a (left) connection
\be 
\na:M\to M\oa\O^1\;: & \qquad & \na(\mu f) = (\na\mu)f + \mu\oa df\;,\\
\na:M^\ast\to\O^1\oa M^\ast\;: & \qquad & \na(fm) = f(\na m) + df\oa m\;,
\ee
so that 
\be
d\langle m,\mu\rangle =\langle\na m,\mu\rangle +\langle m,\na\mu\rangle\;.
\ee
$\na$ can be extended to tensor fields, though the details will not be given
here. Moreover defining 
\be  
\na(\mu\oa\psi) & := & \mu\oa d\psi+(\na\mu)\psi\;,\\
\na(\psi\oa m) & := & d\psi\oa m + (-1)^r\psi(\na m)\;,
\ee
$\na$ is extended to $M\oa\O$, $\O\oa M^\ast$. 
For $M=\O^1$, $M^\ast=\X$ curvature and torsion are defined by
\be
\na^2:\O^1\to\O^1\oa\O^2\;,\qquad
\Theta:\O^1\to\O^2\;:\quad \Theta(\alpha):=d\alpha+\pi(\na\alpha)\;,
\ee
$\pi$ being given by (\ref{pi}). 

It is worth noting that in ordinary DC, there is a basis of 1-forms,
$\{dx^\mu\}$ so that $d(dx^\mu)=0$. This implies that torsion is completely
reducible to the connection. However in 2nd order DC this is not true,
since in general there is {\em no} coordinate system where $d\xi^{\mu\nu}$
vanish. Below we shall see that a covariant basis of $\O^2$ includes 
components of $\Theta$ irreducible to $\na$. 

Further progress follows by extending the $\bu$-product to $\O^2$ and 
$\O^1\oa\O^1$ 
\be
\begin{array}{c}
(\alpha\oa\beta)\bu\gamma := \alpha\oa(\beta\bu\gamma)\\
\alpha\bu(\beta\oa\gamma) := (\alpha\bu\beta)\oa\gamma
\end{array}
\;,\qquad
(\alpha\oa\beta)\bu(\gamma\oa\delta):=(\alpha\bu\gamma)\oa(\beta\bu\delta)\;,
\ee
\be
\omega\bu(gdf\,h):=(gdf\,h)\bu\omega:=g[f,\omega]h\;,\qquad (\alpha\beta)
\bu(\gamma\delta):=(\alpha\bu\gamma)(\beta\bu\delta)\;,
\ee
$\alpha,\,\beta,\,\gamma,\,\delta\in\O^1,\;\omega\in\O^2,\;f,\,g,\,h\in\A$.
Then it is possible to obtain a right $\A$-linear product of 1-forms
through
\be \alpha\circ\beta:=\alpha\beta-\pi(\na\alpha\bu\beta) \ee
since by noncommutativity the product in $\O$ in general is not right
$\A$-linear in both factors $(\alpha f)\beta\neq\alpha\beta f$. Then a 
lengthy calculation gives the identity 
\be
\Theta(\alpha\bu\beta)=\alpha\circ\beta+\beta\circ\alpha - \alpha\bu
\Theta(\beta)-\Theta(\alpha)\bu\beta-\Theta(\alpha)\bu\Theta(\beta)+
\pi B(\alpha,\beta)\;,\label{th}
\ee
where
\be 
B(\alpha,\beta):=\na(\alpha\bu\beta)-\alpha\bu\na\beta-\na\alpha\bu\beta
+(\na\alpha)\bu(\na\beta) \label{bb}
\ee
is a {\em tensor field}. 

Motivated by (\ref{th}) an {\em antisymmetric} and right $\A$-linear
wedge product $\we$ is defined by
\be
\alpha\we\beta:=\alpha\circ\beta+{1\over2}[\Theta(\alpha)\bu\beta+
\alpha\bu\Theta(\beta)-\Theta(\alpha)\bu\Theta(\beta)+\pi B(\alpha,\beta)]
\label{we}
\ee
an indispensable tool for SyM and for obtaining a covariant basis of
$\O^2$ as well.

We first find such a basis of $\O^1$, by noticing that $\O^1\bu\O^1=:
\O^1_2$ is an ideal of $(\O^1,\bu)$. Its annihilator in $\X$ is a submodule 
$\X_1$ of $\X$, that can be shown to consist of all {\em derivations} of
$\A$. Thus, by introducing complementary projections $p_1,\,p_2$ on $\O^1$ and
their duals $p_1^\ast,\,p_2^\ast$ on $\X$ we write
\bez
\O^1=\O^1_1\oplus\O^1_2\;,\quad \X=\X_1\oplus\X_2\;,\qquad
\O^1_i=p_i(\O^1)\;,\quad \X_i=p_i^\ast(\X)\;.
\eez
Writing
\be
p_1(dx^\mu)=dx^\mu + {1\over2}\xi^{\rho\sigma} P^\mu{}_{\rho\sigma}=:
\tilde{d}x^\mu \label{p1o}
\ee
we get
\be
p_2^\ast(\pa_\mu\pa_\nu) = \pa_\mu\pa_\nu - P^\rho{}_{\mu\nu}\pa_\rho
=:\tilde{\pa}_{\mu\nu} 
\ee
and from this, that $P^\rho{}_{\mu\nu}$ are the components of an ordinary,
symmetric (torsionless)	connection, hence $\{ dt,\,\tilde{d}x^\mu,\,
\xi^{\mu\nu}\}$, $\{\pa_t,\,\pa_\mu,\,\tilde{\pa}_{\mu\nu}\}$ are
covariant bases of $\O^1,\,\X$ respectively. Then
\be
df = dt\,\pa_tf + \tilde{d}x^\mu\pa_\mu f + {1\over2}\xi^{\mu\nu}
\tilde{\pa}_{\mu\nu}f\;. 
\ee
Projections $p_i^\ast$ correspond to the introduction of a connection in
StC on manifolds,where it is defined as a mapping (\cite{Em} p.32)
$\Gamma:\A\to \mbox{symmetric bilinear forms on the manifold}$
\bez  \Gamma(f^2) = 2f\Gamma(f) + 2df\oa df\;.  \eez
In fact $\Gamma(f)=p_2(df)$, where we identify $\xi^{\mu\nu}$ with the 
symmetrized $dx^\mu\oa dx^\nu$. Equivalently, $p_1^\ast$ is a mapping 
of 2nd order differential opertors to derivations, its kernel 
$p_2^\ast(\X)$ giving the distribution of horizontal subspaces on the 
manifold (\cite{Me} section 3).

A deeper look at the concept of a connection, presupposes the following 
remarks: In analogy with the fact that $d$, $\na$ follow the {\em grading}
of forms
\bez
\O^1\To{d}\O^2\To{d}\cdots\;,\qquad\qquad
\O^1\oa\A\To{\na}\O^1\oa\O^1\To{\na}\O^1\oa\O^2\To{\na}\cdots
\eez
we would like to have operations $d\bu$, $\na\bu$ say, with {\em similar}
properties, following the {\em filtration} implied by $\bu$
\bez
\O^1\To{d\bu}\O^1_2\To{d\bu}\cdots\;,\qquad\qquad
\O^1\oa\O^1\To{\na\bu}\O^1\oa\O^1_2\To{\na\bu}\cdots\;.
\eez
More precisely, {\em in analogy} to
\bez
d^2=0\;,\qquad\qquad d(f\alpha g)=df\,\alpha g+fd\alpha\,g-f\alpha dg\;,
\eez
\bez 
\na(\alpha\oa\beta) = (\na\alpha)\beta+\alpha\oa d\beta\;,
\eez
we require
\be
d\bu d=0\;,\qquad\qquad	d\bu(f\alpha g)=df\bu\alpha g+fd\bu\alpha g-
f\alpha\bu dg\;, \label{dbu}
\ee
\be
\na\bu(\alpha\oa\beta) = \na\alpha\bu\beta + \alpha\oa d\bu\beta\;.
\label{nbu}
\ee
It turns out that (\ref{dbu}), (\ref{nbu}) specify $d\bu$, $\na\bu$
uniquely and in addition that the {\em curvature-like quantity}
$\na\bu\na\alpha$ is a {\em tensor field}. Similarly,  the {\em 
torsion-like quantity} $(d\bu+\na\bu)\alpha$,
is also a {\em tensor field}, where by definition 
$\na\bu\alpha$ is the $\bu$-product of the two $\oa$-factors in $\na\alpha$.

Using these operations and splitting $\na\alpha$, $\alpha\in\O^1$ as
\be
\na\alpha & = & (p_1\oa p_1)\na p_1(\alpha)+(p_1\oa p_2)\na p_1(\alpha)
+(p_2\oa\mbox{id})\na p_2(\alpha)\nonumber\\
&&+(p_1\oa\mbox{id})\na p_2(\alpha)+(p_2\oa\mbox{id})\na p_1(\alpha)
\label{sc}
\ee
we can show (i) that 
\be 
K(\alpha) := - (p_2\oa\mbox{id})\na p_1(\alpha)\;,\qquad
L(\alpha) := - (p_1\oa\mbox{id})\na p_2(\alpha)
\ee
are {\em tensor fields}, (ii) the 2nd and 3rd term contain respectively
the curvature-like tensor field
\be
A(\alpha):=-{1\over2}(p_1\oa\mbox{id})\na\bu\na p_1(\alpha)
\ee
and $B$ in (\ref{bb}), (iii) the remaining parts can be expressed via
$K,\,L,\,A,\,B$ and the 1st term in (\ref{sc}) which represents a
{\em 1st order connection}. More explicitely, we write
\be
\lefteqn{\na(\tilde{d}x^\mu) =}&&\hspace{15.5cm} \nonumber\\
&&\;\;\; -\tilde{d}x^\nu\oa\tilde{d}x^\rho
\Gamma^\mu{}_{\rho\nu}-{1\over2}\tilde{d}x^\nu\oa\xi^{\rho\sigma}
\Gamma^\mu{}_{\widetilde{\rho\sigma}\,\nu}-{1\over2}\xi^{\rho\sigma}
\oa\tilde{d}x^\nu\Gamma^\mu{}_{\nu\,\widetilde{\rho\sigma}}-
{1\over4}\xi^{\rho\sigma}\oa\xi^{\kappa\lambda}\Gamma^\mu{}
_{\widetilde{\kappa\lambda}\,\widetilde{\rho\sigma}}\;,
\ee
\be
\lefteqn{\na(\xi^{\mu\nu}) =}&& \hspace{15.5cm}\nonumber\\
&& -\tilde{d}x^\rho\oa\tilde{d}x^\sigma\Gamma^{\widetilde{\mu\nu}}{}
_{\sigma\rho}-{1\over2}\tilde{d}x^\rho\oa\xi^{\kappa\lambda}\Gamma
^{\widetilde{\mu\nu}}{}_{\widetilde{\kappa\lambda}\,\rho}
-{1\over2}\xi^{\rho\sigma}\oa\tilde{d}x^\kappa\Gamma
^{\widetilde{\mu\nu}}{}_{\kappa\,\widetilde{\rho\sigma}}-{1\over4}
\xi^{\rho\sigma}\oa\xi^{\kappa\lambda}\Gamma^{\widetilde{\mu\nu}}
{}_{\widetilde{\kappa\lambda}\,\widetilde{\rho\sigma}}\;,
\ee
where $t$ is considered as the $(N+1)$-coordinate and an index of the
form $\widetilde{\rho\sigma}$ is symmetric and corresponds to a 
$\xi^{\rho\sigma}$-component. Then lengthy calculations show that all
$\Gamma$'s are expressed via $\Gamma^\mu{}_{\nu\rho}$ and the components
of $A,\,B,\,K,\,L$. The expressions are complicated and will be given
elsewhere. 

Finally the torsion-like field
\be     S(\alpha):=(d\bu+\na\bu)p_1(\alpha)   \ee
has components (cf.\ (\ref{p1o}))
\be 
S^\mu{}_{\widetilde{\alpha\beta}} = -{1\over2}\Gamma^\mu{}_{(\alpha\beta)}
+P^\mu{}_{\alpha\beta}\;.
\ee
Thus we conclude: (a) {\em pure connection is ordinary connection}
$\Gamma^\alpha{}_{\beta\gamma}$ only, all corrections due to 2nd order DC 
being the tensor fields $K,\,L,\,B,\,A$, (b) any connection induces a 
{\em projection} on $\O^1$, by taking $S=0$, (c) in view of (a), (b), and 
making a {\em minimal} deformation of ordinary differential geometry, we make  
the {\em minimal choice}
\be     K=L=B=A=S=0\;.  \ee
To simplify the calculations we further assume that $\Gamma^\alpha{}_{\beta
\gamma}=\Gamma^\alpha{}_{\gamma\beta}$. 

To get a covariant basis of $\O^2$, we notice that
\bez
\{\tilde{d}x^\mu\circ\tilde{d}x^\nu,\; \tilde{d}x^\mu\circ\xi^{\rho\sigma},
\;\xi^{\mu\nu}\circ\tilde{d}x^\rho,\; \xi^{\mu\nu}\circ\xi^{\rho\sigma} \}
\eez
together with $dt\tilde{d}x^\mu,\,dt\xi^{\mu\nu}$ form a covariant set 
spanning $\O^2$. However it is not a linearly independent set,
because by applying $d$ to (\ref{bu}) it can be shown that the following
relations hold	$d\xi^{\mu\nu}=[dx^\mu,dx^\nu]_+$
\be
[dx^{(\mu},\xi^{\nu\rho)}]_+ =0\;,\qquad
[\xi^{(\mu\nu},\xi^{\rho\sigma)}]_+=0\;,
\ee
\be
dtdt=0\;,\qquad [dt,dx^\mu]_+=0\;,\qquad [dt,\xi^{\mu\nu}]_+=0\;,
\ee
where $[\,,]_+$ denotes the anticommutator with respect to the graded 
product. Moreover no other 2-form relations exist. Putting 
\bez 
\Theta^\mu:=\Theta(\tilde{d}x^\mu)\;,\qquad
\Theta^{\mu\nu}:=\Theta(\xi^{\mu\nu})
\eez
we can show that the {\em torsion-components irreducible to the
connection} are
\bez
\Theta^{\mu\nu}\;,\qquad \Theta^{\mu[\nu}\bu\Theta^{\rho]\sigma}\;,
\qquad \Theta^{\mu[\nu}\bu\tilde{d}x^{\rho]}\;.
\eez  
These, together with $\tilde{d}x^\mu\we\tilde{d}x^\nu$, $\tilde{d}x^\mu
\we\xi^{\rho\sigma}$, $\xi^{\mu\nu}\we\xi^{\rho\sigma}$, $dt\we\tilde{d}x^
\mu$, $dt\we\xi^{\mu\nu}$ form a covariant basis of $\O^2$ via the
relations implied by (\ref{we})
\be
\tilde{d}x^\mu\we\tilde{d}x^\nu & = & \tilde{d}x^\mu\circ\tilde{d}x^\nu
-{1\over2}\Theta^{\mu\nu} +{1\over2}\Theta^{(\mu}\bu\tilde{d}x^{\nu)}\;,\\
\tilde{d}x^\mu\we\xi^{\rho\sigma} & = & \tilde{d}x^\mu\circ\xi^{\rho\sigma}
+{1\over6}(\Theta^{\rho[\sigma}\bu\tilde{d}x^{\mu]}+\Theta^{\sigma[\rho}
\bu\tilde{d}x^{\mu]})\;,\\
\xi^{\mu\nu}\we\xi^{\rho\sigma} & = & \xi^{\mu\nu}\circ\xi^{\rho\sigma}+
{1\over6}(\Theta^{\mu[\nu}\bu\Theta^{\rho]\sigma}+
\Theta^{\mu[\nu}\bu\Theta^{\sigma]\rho})
\ee
where we can find that
\be
\Theta^\mu\bu\tilde{d}x^\nu={1\over2}(\xi^{\rho\nu}\we\xi^{\alpha\beta}+
{1\over12}\Theta^{\nu[\rho}\bu\Theta^{\beta]\alpha})R^\mu{}_{\alpha\beta
\rho}\;,
\ee
$R^\mu{}_{\alpha\beta\rho}$ being the Riemann tensor of $\Gamma^\alpha{}_
{\beta\gamma}$.
\vspace*{1cm}

\centerline{\bf 4. SYMPLECTIC MECHANICS}
\renewcommand{\theequation} {4.\arabic{equation}}
\setcounter{equation}{0}
\vspace*{.5cm}

\noindent
Motivated by the discussion in section 1 (cf.\ (\ref{df}), (\ref{ito})),
we consider the special case $N=2n$ and
\be \xi^{\mu\nu} = -dt b^{\mu\nu}\;, \label{xb}\ee
$b^{\mu\nu}$ being a symmetric matrix function (under (\ref{xb}), (\ref{df})
reduces formally to (\ref{ito})). Application of the results of section 3
gives 
\be 
\Theta^{\mu\nu}\bu\Theta^{\rho\sigma}=
\Theta^{\mu\nu}\bu\tilde{d}x^\rho=0\;,
\ee
\be
\Theta^\mu={1\over2}dx^\rho dt b^{\alpha\beta} R^\mu{}_{\alpha\beta\rho}\;,
\qquad
\Theta^{\mu\nu}=dt dx^\kappa \na_\kappa b^{\mu\nu}
\ee
so that 
\be 
\tilde{d}x^\mu\we dt & = & dx^\mu dt\;, \\
\tilde{d}x^\mu\we\tilde{d}x^\nu & = & \tilde{d}x^\mu\tilde{d}x^\nu
+dtdx^\kappa(b^{\nu\alpha}\Gamma^\mu{}_{\kappa\alpha}-
{1\over2}\na_\kappa b^{\mu\nu})\label{4.5}
\ee
form a basis of $\O^2$ --- here $\na_\kappa$ is defined by 
$\Gamma^\alpha{}_{\beta\gamma}$. Moreover, $\we$ can be extended to any
forms in $\O$ (\cite{DiTz} eq.(4.25)). Then writing for any $\omega\in
\O^2$
\be
\omega={1\over2}\tilde{d}x^\mu\we\tilde{d}x^\nu \omega_{\mu\nu}
+dtdx^\mu\omega_\mu\label{4.6}
\ee
a lengthy calculation gives that {\em $d\omega=0$ is equivalent to
$\omega_{\mu\nu}$ being closed in the ordinary sense and}
\be
\omega_\mu=-{1\over2}b^{\nu\rho}\na_\nu\omega_{\rho\mu} +\pa_\mu H\label{4.7}
\ee
for some function $H$. Actually, (\ref{4.5}), (\ref{4.7}) imply that
(\ref{4.6}) reduces to (\ref{1.4}) as in conventional SyM. With $\omega$
having maximal rank, and applying the Hamiltonian scheme (\ref{1.5}),
(\ref{hd}) we get Hamilton's equations in the 2nd order DC 
$X_H(A)=0\iff$
\be
\pa_t A = -(\{H,A\}+F_\mu\omega^{\mu\nu}\pa_\nu A)+{1\over2}\pa_\mu
(b^{\mu\nu}\pa_\nu A)+{1\over2}b^{\mu\nu}\Gamma^\rho{}_{\rho\nu}
\pa_\mu A\;,\label{4.8}
\ee
\be
F_\mu := -{1\over2}\na_\nu(b^{\nu\rho}\omega_{\rho\mu})\;,\qquad
\omega^{\mu\rho}\omega_{\nu\rho} =\delta^\mu_\nu\;. 
\ee
Clearly (\ref{4.8}) is identical to (\ref{fp}), provided
\be 
F_\mu=\pa_\mu F\;,\qquad H_0 = H+F\;,\qquad \Gamma^\rho{}_{\rho\mu}=
-\beta\pa_\mu H_0\label{4.10}
\ee
for some function $F$. It can be shown that {\em (\ref{4.10}) are equivalent
to $\omega_{\mu\nu}$ being harmonic with respect to the generalized
Laplace-Beltrami operator of $b^{\mu\nu}$\footnote{If $b^{\mu\nu}$ is 
nondegenerate, and $\Gamma^\rho{}_{\mu\nu}$ its metric connection then
$\omega_{\mu\nu}$ must be harmonic.} and to the canonical volume 
$e^{-\beta H}dx^1\cdots dx^{2n}$ being covariantly constant}\footnote{In
statistical mechanics this is the Gibbs measure given by the 
Maxwell-Boltzmann probability density $e^{-\beta H}$, $H$ being the system's
Hamiltonian.} (see \cite{DiTz} section 4).
\vspace*{1cm}

\centerline{\bf 5. DISCUSSION}
\renewcommand{\theequation} {5.\arabic{equation}}
\setcounter{equation}{0}
\vspace*{.5cm}

\noindent
Below we comment on a possible realization of 2nd order DC: By (\ref{2.2})
for $f_1,\ldots,f_k\in\A$
\be 
(d_uf_1\bu_u\cdots\bu_ud_uf_k)(a,b)=\Delta f_1(a,b)\cdots\Delta f_k(a,b)\;,
\qquad \Delta f_i(a,b):=f_i(b)-f_i(a)\;.\label{5.1}
\ee
In general (\ref{5.1}) cannot be valid in a DC $(\O,d,\bu)$, being incompatible
with its defining relations, (equivalently with the projection
$(\O_u,d_u,\bu_u)\to(\O,d,\bu)$, cf.\ (\ref{2.1}), (\ref{pi})). However
if it is to be somehow retained, it leads to some interesting suggestions:

For instance in the ordinary DC the {\em nontrivial} defining relation is
 $df\bu dg=0$, so that in one dimension 
\be  (dx\bu dx)(a,b)=0 \quad \Longrightarrow \quad (b-a)^2=0\;. \ee
Excluding the {\em trivial} case $b \ = \ a$,
we may reinterpret the above result by saying that $a,\,b$ {\em differ by
an infinitesimal of first order}. When $a,\,b$ are not neighbouring, 
for $\alpha\in\O^1$, $\alpha(a,b)$  may be interpreted as a kind of integral
\bez 
\int^b_a df\bu dg = \lim\sum_{k=1}^n\Delta f(x_{k-1},x_k)
\Delta g(x_{k-1},x_k) 
\eez
the limit being taken with respect to the width of the partition
$a=x_0 < x_1 < \cdots < x_n=b$. Thus $df\bu df=0$ means, that $f$
{\em must be of zero quadratic variation, hence of bounded variation}. Since
$f,\,g$ are continuous, this ensures the existence of $\int^b_a f dg$ as 
a Riemann-Stieltjes integral.    

In the case of the 2nd order DC, (\ref{bu}) gives similarly $(b-a)^3=0$,
or by the above reasoning, {\em that $a,\,b$ differ by an infinitesimal of 
2nd order}. For arbitrary $a,\,b$
\bez
\int^b_a df\bu dg\bu dh=0\;\Longrightarrow\;\lim\sum_{k=1}^n
\Delta f(x_{k-1},x_k)\Delta g(x_{k-1},x_k)\Delta h(x_{k-1},x_k)=0\;,
\eez
hence $f,\,g,\,h$ must have zero {\em cubic variation}, hence
{\em finite quadratic variation} given by $\int^b_a df\bu df$ etc.
This corresponds directly to stochastic integration, (see (\ref{sm}),
(\ref{1.3}) and the comments following (\ref{2.6})). 

The above {\em nonrigorous} remarks suggest that StC may be a possible
realization of 2nd DC, but it is a still unsolved problem, whether
it is the only one, and much remains to be done in this direction.

\renewcommand{\refname}{

\end{document}